\documentclass[a4paper,12pt]{article}
\begin{document}

\title{Manifolds over Cayley-Dickson algebras and their immersions.}
\author{Ludkovsky S.V.}
\date{15 March 2012}
\maketitle

\begin{abstract}
Holomorphic manifolds over Cayley-Dickson algebras are defined and
their embeddings and immersions are studied.
\end{abstract}

\section{Introduction.}
Real and complex manifolds are widely used in different branches of
mathematics \cite{henlei,kodai,michor}. On the other hand,
Cayley-Dickson algebras ${\cal A}_r$, particularly, the quaternion
skew field ${\bf H}={\cal A}_2$ and the octonion algebra ${\bf
O}={\cal A}_3$, have found many-sided applications not only in
mathematics, but also in theoretical physics
\cite{brdeso,gilmurr,girard,guespr,harvey,lawmich,oystaey}.
Functions of Cayley-Dickson variables were studied earlier
\cite{ludoyst,ludfov,lufscdvm,lufjmsrf}. Their
super-differentiability was defined in terms of representing them
words and phrases as a differentiation which is real-linear,
additive and satisfying Leibniz' rule on an algebra of phrases over
${\cal A}_r$. A super-differentiable function on a domain $U$ in
${\cal A}_r^n$ or $l_2({\cal A}_r)$ of ${\cal A}_r$-variables is
also called holomorphic. \par This article is devoted to
investigations of ${\cal A}_r$-holomorphic manifolds. Their
embeddings and immersions are studied. Results and notations of
previous papers \cite{ludoyst,ludfov,lufscdvm,lufjmsrf} are used
below.
\par Main results of this paper are obtained for the first time.

\section{Manifolds over Cayley-Dickson algebras}
{\bf 1. Definitions and Notes.} An ${\bf R}$ linear space $X$ which
is also left and right ${\cal A}_r$ module will be called an ${\cal
A}_r$ vector space. We present $X$ as the direct sum
\par $(DS)$ $\quad X=X_0i_0\oplus ... \oplus X_m i_m\oplus ...$, where
$X_0$,...,$X_m,...$ are pairwise isomorphic real linear spaces,
where $i_0,...,i_{2^r-1}$ are generators of the Cayley-Dickson
algebra ${\cal A}_r$ such that $i_0=1$, $i_k^2=-1$ and
$i_ki_j=-i_ji_k$ for each $k\ge 1$ and $j\ge 1$ so that $k\ne j$,
$~2\le r$.
\par  Let $X$ and $Y$ be two
$\bf R$ linear normed spaces which are also left and right ${\cal
A}_r$ modules, where $1\le r$, such that  \par $(1)$ $0\le \| ax
\|_X \le |a| \| x \|_X $ and $ \| x a \|_X\le |a| \| x \|_X $ and
\par $(2)$ $\| ax_j \|_X =
|a| \| x_j \|_X$ and \par $(3)$ $\| x+y \|_X \le \| x \|_X + \| y \|_X $ \\
for all $x, y\in X$ and $a\in {\cal A}_r$ and $x_j\in X_j$. Such
spaces $X$ and $Y$ will be called ${\cal A}_r$ normed spaces.
\par Suppose that $X$ and $Y$ are two normed spaces over the Cayley-Dickson algebra ${\cal
A}_v$. A continuous $\bf R$ linear mapping $\theta : X\to Y$ is
called an $\bf R$ linear homomorphism. If in addition $\theta
(bx)=b\theta (x)$ and $\theta (xb)=\theta (x)b$ for each $b \in
{\cal A}_v$ and $x\in X$, then $\theta $ is called a homomorphism of
${\cal A}_v$ (two sided) modules $X$ and $Y$.
\par If a homomorphism is injective, then it is called an embedding
($\bf R$ linear or for ${\cal A}_v$ modules correspondingly).
\par If a homomorphism $h$ is bijective and from $X$ onto $Y$ so that its inverse
mapping $h^{-1}$ is also continuous, then it is called an
isomorphism ($\bf R$ linear or of ${\cal A}_v$ modules
respectively).

\par {\bf 2. Definitions.} We say that a real vector space $Z$ is supplied with a scalar
product if a bi-$\bf R$-linear bi-additive mapping $<,>: Z^2\to {\bf
R}$ is given satisfying the conditions:
\par $(1)$ $<x,x> ~ \ge 0$, $ ~ <x,x>=0$ if and only if $x=0$;
\par $(2)$ $<x,y>=<y,x>$;
\par $(3)$ $<ax+by,z>=a<x,z>+b<y,z>$ for each real numbers $a, b\in {\bf R}$ and vectors $x, y,
z\in Z$.
\par Then an ${\cal A}_r$ vector space $X$ is supplied with an
${\cal A}_r$ valued scalar product, if a bi-${\bf R}$-linear
bi-${\cal A}_r$-additive mapping $<*,*>: X^2\to {\cal A}_r$ is given
such that
\par $(4)$ $\quad <f,g> = \sum_{j,k} <f_j,g_k>i_j^*i_k$, \\
where $f=f_0i_0+...+f_mi_m+...$, $ ~f, g\in X$, $ ~ f_j, g_j \in
X_j$, each $X_j$ is a real linear space with a real valued scalar
product, $(X_j, <*,*>)$ is real linear isomorphic with $(X_k,
<*,*>)$ and $<f_j,g_k>\in {\bf R}$ for each $j, k$. The scalar
product induces the norm: \par $(5)$ $\| f \| := \sqrt{<f,f>}$.
\par An ${\cal A}_r$ normed space or an
${\cal A}_r$ vector space with ${\cal A}_r$ scalar product complete
relative to its norm will be called an ${\cal A}_r$ Banach space or
an ${\cal A}_r$ Hilbert space respectively.
\par A Hilbert space $X$ over ${\cal A}_r$ is denoted by $l_2(\lambda ,{\cal
A}_r)$, where $\lambda $ is a set of the cardinality $card (\lambda
)\ge \aleph _0$ which is the topological weight of $X_0$, i.e.
$X_0=l_2(\lambda ,{\bf R})$. \par A mapping $f: U\to l_2(\lambda
,{\cal A}_r)$ can be written in the form $$f(z) = \sum_{j\in \lambda
} f^j(z)e_j,$$ where $\{ e_j: j \in \lambda  \} $ is an orthonormal
basis in the Hilbert space $l_2(\lambda ,{\cal A}_r)$, $~U$ is a
domain in $l_2(\psi ,{\cal A}_r)$, $~f^j(z)\in {\cal A}_r$ for each
$z\in U$ and every $j\in \lambda $. If $f$ is Frech\'et
differentiable over $\bf R$ and each function $f^j(z)$ is
holomorphic by each Cayley-Dickson variable $\mbox{}_kz$ on $U$,
then $f$ is called holomorphic on $U$, where $$z = \sum_{k\in \psi }
\mbox{}_kzq_k,$$ while $\{ q_k: ~ k\in \psi \} $ denotes the
standard orthonormal basis in $l_2(\psi ,{\cal A}_r)$, $\mbox{
}_kz\in {\cal A}_r$.

\par {\bf 3. Definition.} Let $M$ be a set such that $M=\bigcup_jU_j$,
$M$ is a Hausdorff topological space, each $U_j$ is open in $M$,
$\phi _j: U_j\to \phi _j(U_j)\subset X$ are homeomorphisms, $\phi
_j(U_j)$ is open in ${\cal A}_r$ for each $j$, if $U_i\cap U_j\ne
\emptyset $, the transition mapping $\phi _i\circ \phi _j^{-1}$ is
${\cal A}_r$-holomorphic on its domain, where $X$ is either ${\cal
A}_r^m$ with $m\in {\bf N}$ or a Hilbert space $l_2(\lambda ,{\cal
A}_r)$ over the Cayley-Dickson algebra ${\cal A}_r$. Then $M$ is
called the ${\cal A}_r$-holomorphic manifold.

\par {\bf 4. Proposition.} {\it Let $M$ be an ${\cal A}_r$ holomorphic
manifold. Then there exists a tangent bundle $TM$ which has the
structure of an ${\cal A}_r$ holomorphic manifold such that each
fibre $T_xM$ is the vector space over the Cayley-Dickson algebra
${\cal A}_r$.}
\par {\bf Proof.} The Cayley-Dickson algebra ${\cal A}_r$ has the real shadow,
which is the Euclidean space $\bf R^{2^r}$, since ${\cal A}_r$ is
the algebra over $\bf R$. Therefore, a manifold $M$ has also a real
manifold structure. Each ${\cal A}_r$ holomorphic mapping is
infinite differentiable in accordance with Theorems 2.15 and 3.10 in
\cite{ludoyst,ludfov}. Then the tangent bundle $TM$ exists, which is
$C^{\infty }$-manifold such that each fibre $T_xM$ is a tangent
space, where $x\in M$, $~T$ is the tangent functor. If $At (M)= \{
(U_j,\phi _j): j \} $, then $At (TM)= \{ (TU_j, T\phi _j): j \} $,
$TU_j=U_j\times X$, where $X$ is the ${\cal A}_r$ vector space on
which $M$ is modeled, $T(\phi _j\circ \phi _k^{-1})=(\phi _j\circ
\phi _k^{-1}, D(\phi _j\circ \phi _k^{-1}))$ for each $U_j\cap
U_k\ne \emptyset $. Each transition mapping $\phi _j\circ \phi
_k^{-1}$ is ${\cal A}_r$ holomorphic on its domain, then its
(strong) differential coincides with the super-differential $D(\phi
_j\circ \phi _k^{-1})= D_z(\phi _j\circ \phi _k^{-1})$, since
${\tilde \partial } (\phi _j\circ \phi _k^{-1})=0$. Therefore, the
super-differential $D(\phi _j\circ \phi _k^{-1})$ is $\bf R$-linear
and ${\cal A}_r$-additive, hence it is the automorphism of the
${\cal A}_r$ vector space $X$. But $D_z(\phi _j\circ \phi _k^{-1})$
is ${\cal A}_r$ holomorphic as well, consequently, $TM$ is the
${\cal A}_r$ holomorphic manifold.

\par {\bf 5. Definitions.} A $C^1$-mapping $f: M\to N$ is called
an immersion, if $rang (df|_x: T_xM\to T_{f(x)}N) = m_M$ for each
$x\in M$, where $m_M := dim_{\bf R}M$. An immersion $f: M\to N$ is
called an embedding, if $f$ is bijective.

\par {\bf 6. Theorem.} {\it Let $M$ be a compact
${\cal A}_r$ holomorphic manifold, $dim_{{\cal A}_r}M = m<\infty $,
where $2\le r$. Then there exists an ${\cal A}_r$ holomorphic
embedding $\tau : M\hookrightarrow {\cal A}_r^{2m+1}$ and an ${\cal
A}_r$ holomorphic immersion $\theta : M\to {\cal A}_r^{2m}$
correspondingly. Each continuous mapping $f: M\to {\cal A}_r^{2m+1}$
or $f: M\to {\cal A}_r^{2m}$ can be approximated by $\tau $ or
$\theta $ relative to the norm $\| * \|_{C^0}$. If $M$ is a
paracompact ${\cal A}_r$ holomorphic manifold with countable atlas
on $l_2(\lambda ,{\cal A}_r)$, where $card (\lambda )\ge \aleph _0$,
then there exists a holomorphic embedding $\tau : M\hookrightarrow
l_2(\lambda ,{\cal A}_r)$.}
\par {\bf Proof.} Let at first $M$ be compact.
Since $M$ is compact, then it is finite dimensional over ${\cal
A}_r$, $dim_{{\cal A}_r}M = m\in \bf N$, such that $dim_{{\cal
A}_r}M = 2^rm$ is its real dimension. Take an atlas $At' (M)$
refining the initial atlas $At (M)$ of $M$ such that $({U'}_j, \phi
_j)$ are charts of $M$, where each ${U'}_j$ is ${\cal A}_r$
holomorphic diffeomorphic to an interior of the unit ball $Int
(B({\cal A}_r^m,0,1))$, where $B({\cal A}_r^m,y,\rho ) := \{ z\in
{\cal A}_r^m: |z-y|\le \rho \} $. In view of compactness of the
manifold $M$ a covering $\{ {U'}_j: j \} $ has a finite subcovering,
hence $At' (M)$ can be chosen finite. Denote for convenience the
latter atlas as $At (M)$. Let $(U_j, \phi _j)$ be the chart of the
atlas $At (M)$, where $U_j$ is open in $M$, hence $M\setminus U_j$
is closed in $M$.
\par  Consider the space ${\cal A}_r^m\times {\bf R}$
as the $\bf R$-linear space $\bf R^{2^rm+1}$. The unit sphere
$S^{2^rm}:=S ({\bf R}^{2^rm+1},0,1) := \{ z\in {\bf R}^{2^rm+1}:$
$|z|=1 \} $ in ${\cal A}_r^m\times \bf R$ can be supplied with two
charts $(V_1, \phi _1)$ and $(V_2, \phi _2)$ such that $V_1:=
S^{2^rm}\setminus \{ 0,...,0, 1 \} $ and $V_2:=S^{2^rm}\setminus \{
0,...,0, - 1 \} $, where $\phi _1$ and $\phi _2$ are stereographic
projections from poles $\{ 0,...,0, 1 \} $ and $ \{ 0,...,0, -1 \} $
of $V_1$ and $V_2$ respectively onto ${\cal A}_r^m$. The conjugation
is given by the formula: $$z^* = - (2^r-2)^{-1}
\sum_{p=0}^{2^r-1}(i_pz)i_p $$ in ${\cal A}_r^m$, which provides
$z^*$ in the $z$-representation. Therefore $\phi _1\circ \phi
_2^{-1}$ written in the $z$-representation is the ${\cal A}_r$
holomorphic diffeomorphism in ${\cal A}_r^m\setminus \{ 0 \} $, i.e.
the super-differential $D_z (\phi _1\circ \phi _2^{-1})$ exists,
where $i_0,...,i_{2^r-1}$ are the standard generators of ${\cal
A}_r$. Thus the unit sphere $S^{2^rm}$ can be supplied with the
structure of the ${\cal A}_r$ holomorphic manifold.
\par Therefore, there exists an ${\cal A}_r$ holomorphic mapping
$\psi _j$. That is locally $z$-analytic of $M$ into the unit sphere
$S^{2^rm}$ such that $\psi _j: U_j\to \psi _j(U_j)$ is the ${\cal
A}_r$ holomorphic diffeomorphism onto the subset $\psi _j(U_j)$ in
$S^{2^rm}$. Using of such mappings $\psi _j$ is sufficient, where
$\psi _j$ can be considered as components of a holomorphic
diffeomorphism: $\psi : M\to (S^{2^rm})^n$ with $n$ equal to the
number of charts.  There is an embedding of ${\cal A}_r^m\times \bf
R$ into ${\cal A}_r^{m+1}$. Then the mapping $\psi (z):=(\psi
_1(z),...,\psi _n(z))$ is the embedding into $(S^{2^rm})^n$ and
hence into ${\bf K}^{n{m+1}}$, since the rank is $rank [d_z\psi
(z)]=2^rm$ at each point $z\in M$. Indeed, the rank is $rank
[d_z\psi _j(z)]=2^rm$ for each $z\in U_j$ and the dimension is
bounded from above $dim_{{\cal A}_r}\psi (U_j)\le dim_{{\cal A}_r} M
= m$. Moreover, $\psi (z)\ne \psi (y)$ for each $z\ne y\in U_j$,
since $\psi _j(z)\ne \psi _j(y)$. If $z\in U_j$ and $y\in M\setminus
U_j$, then there exists a number $l\ne j$ so that $y\in U_l\setminus
U_j$, $\psi _j(z)\ne \psi _j(y)=x_j$.
\par Let $M\hookrightarrow {\cal A}_r^N$ be the ${\cal A}_r$
holomorphic embedding as above. There is also the ${\cal A}_r$
holomorphic embedding of $M$ into $(S^{2^rm})^n$ as it is shown
above, where $(S^{2^rm})^n$ is the ${\cal A}_r$ holomorphic manifold
as the product of ${\cal A}_r$ holomorphic manifolds. \par Let
$P{\bf R}^n$ denote the real projective space formed from the
Euclidean space ${\bf R}^{n+1}$, denote by $\phi : {\bf R}^{n+1}\to
P{\bf R}^n$ the corresponding projective mapping. Geometrically
$P{\bf R}^n$ is considered as $S^n/\tau $, where $S^n:= \{ y\in {\bf
R}^{n+1}: ~ \| y \| =1 \} $ is the unit sphere in ${\bf R}^{n+1}$,
while $\tau $ is the equivalence relation making identical two
spherically symmetric points, i.e. points belonging to the same
straight line containing zero and intersecting the unit sphere.
Consider ${\cal A}_r^n$ as the algebra of all $n\times n$ diagonal
matrices $A = diag (a_1,...,a_n)$ with entries $a_1,...,a_n\in {\cal
A}_r$. Then ${\cal A}_r^n$ is isomorphic with the tensor product of
algebras ${\cal A}_r^n = {\cal A}_r \otimes_{\bf R} {\bf R}^n$ over
the real field, where ${\bf R}^n$ is considered as the algebra of
all diagonal $n\times n$ matrices $C=diag (b_1,..,b_n)$ with entries
$b_1,...,b_n\in {\bf R}$. Put $P{\cal A}_r^n=\phi ({\cal
A}_r\otimes_{\bf R} {\bf R}^{n+1})$ to be the Cayley-Dickson
projective space, where $\phi $ is extended from ${\bf R}^{n+1}$
onto ${\cal A}_r \otimes_{\bf R} {\bf R}^{n+1}$ so that $\phi
(ax)=a\phi (x)$ and $\phi (xa)=\phi (x)a$ for each $a\in {\cal A}_r$
with $|a|=1$ and every $x\in {\bf R}^{n+1}$, also $\phi
(x_0i_0+...+x_{2^r-1}i_{2^r-1})=\phi (x_0)i_0\alpha _0+...+\phi
(x_{2^r-1})i_{2^r-1}\alpha _{2^r-1}$ for each non-zero vector $
x=x_0i_0+...+x_{2^r-1}i_{2^r-1}\in {\cal A}_r^{n+1}$, where $\alpha
_j := \| x_j \| / \| x \| $, $x_j\in {\bf R}^{n+1}$ for each $j$,
$~\| x \| ^2 = \| x_0 \| ^2 +... + \| x_{2^r-1} \| ^2$.
\par If $z\in P{\cal A}_r^n$, then by our definition $\phi ^{-1}(z)$
is the ${\cal A}_r$ straight line in ${\cal A}_r^{n+1}$. To each
element $x\in {\cal A}_r^{n+1}$ we pose an ${\cal A}_r$ straight
line $< {\cal A}_r,x \} := \phi ^{-1}(\phi (x))$. That is
 the bundle of all ${\cal
A}_r$ straight lines $< {\cal A}_r,x \}$ in ${\cal A}_r^N$ is
considered, where $x\in {\cal A}_r^N$, $x\ne 0$, so that $<{\cal
A}_r,x \} $ is the ${\cal A}_r$ vector space of dimension $1$ over
${\cal A}_r$, which has the real shadow isomorphic with ${\bf
R}^{2^r}$.
\par Fix the standard orthonormal base $ \{ e_1,...,e_N \} $ in
${\cal A}_r^N$ and projections on ${\cal A}_r$-vector subspaces
relative to this base
$$P^L(x):=\sum_{e_j\in L}x_je_j$$ for the ${\cal A}_r$ vector span
$L=span_{{\cal A}_r} \{ e_i:$ $i\in \Lambda _L \} $, $ ~ \Lambda
_L\subset \{ 1,...,N \} $, where $$x=\sum_{j=1}^Nx_je_j,$$  $x_j\in
{\cal A}_r$ for each $j$, $ ~ e_j=(0,...,0,1,0,...,0)$ with $1$ at
$j$-th place. In this base consider the ${\cal A}_r$-Hermitian
scalar product $$<x,y> := \sum_{j=1}^Nx_j^*y_j.$$ Let $l\in {\cal
A}_rP^{N-1}$, take an ${\cal A}_r$-hyperplane denoted by $({\cal
A}_r^{N-1})_l$ and given by the condition: $<x,y>=0$ for each $x\in
({\cal A}_r^{N-1})_l$ and $y\in l$. Take a vector $0\ne [l]\in {\cal
A}_r^N$ as a representative characterizes the equivalence class $l=<
{\cal A}_r, [l] \} $ of unit norm $\| [l] \| =1$. Then the
orthonormal base $\{ q_1,...,q_{N-1} \} $ in $({\cal A}_r^{N-1})_l$
and the vector $[l]=:q_N$ compose the orthonormal base $\{
q_1,...,q_N \} $ in ${\cal A}_r^N$. This provides the ${\cal A}_r$
holomorphic projection $\pi _l: {\cal A}_r^N\to ({\cal
A}_r^{N-1})_l$ relative to the orthonormal base $ \{ q_1,...,q_N \}
$. Indeed, the operator $\pi _l$ is ${\cal A}_r$ left $\pi
_l(bx_0)=b\pi _l(x_0)$ and also right $\pi _l(x_0b)=\pi _l(x_0)b$
linear for each $x_0\in X_0$ and $b\in {\cal A}_r$, but certainly
non-linear relative to ${\cal A}_r$. Therefore the mapping $\pi _l$
is ${\cal A}_r$ holomorphic.
\par To construct an immersion it is sufficient, that each
projection $\pi _l: T_xM\to ({\cal A}_r^{N-1})_l$ has $ker [d(\pi
_l(x))]= \{ 0 \} $ for each $x\in M$. The set of all points $x\in M$
for which $ker [d(\pi _l(x))] \ne \{ 0 \} $ is called the set of
forbidden directions of the first kind. Forbidden are those and only
those directions $l\in {\cal A}_rP^{N-1}$ for which there exists a
point $x\in M$ such that $l'\subset T_xM$, where $l'=[l]+z$, $z\in
{\cal A}_r^N$. The set of all forbidden directions of the first kind
forms the ${\cal A}_r$ holomorphic manifold $Q$ of ${\cal A}_r$
dimension $(2m-1)$ with points $(x,l)$, where $x\in M$, $ ~ l\in
{\cal A}_rP^{N-1}$, $ ~ [l]\in T_xM$. Take the mapping $g: Q\to
{\cal A}_rP^{N-1}$ given by $g(x,l):=l$. Then this mapping $g$ is
${\cal A}_r$ holomorphic.
\par Each paracompact manifold $A$ modeled
on ${\cal A}_r^p$ can be supplied with the Riemann manifold
structure also. Therefore, on a manifold $A$ there exists a Riemann
volume element. In view of the Morse theorem $\mu (g(Q))=0$, if
$N-1>2m-1$, that is, $2m<N$, where $\mu $ is the Riemann volume
element in ${\cal A}_rP^{N-1}$. In particular, $g(Q)$ is not
contained in ${\cal A}_rP^{N-1}$ and there exists $l_0\notin g(Q)$,
consequently, there exists $\pi _{l_0}: M\to ({\cal
A}_r^{N-1})_{l_0}$. This procedure can be prolonged, when $2m<N-k$,
where $k$ is the number of the step of projection. Hence $M$ can be
immersed into ${\cal A}_r^{2m}$.
\par Consider now the forbidden directions of the second type:
$l\in {\cal A}_rP^{N-1}$, for which there exist $x\ne y\in M$
simultaneously belonging to $l$ after suitable parallel translation
$[l]\mapsto [l]+z$, $z\in {\cal A}_r^N$. The set of the forbidden
directions of the second type forms the manifold $\Phi
:=M^2\setminus \Delta $, where $\Delta := \{ (x,x):$ $x\in M \} $.
Consider $\psi : \Phi \to {\cal A}_rP^{N-1}$, where $\psi (x,y)$ is
the straight ${\cal A}_r$-line with the direction vector $[x,y]$ in
the orthonormal base. Then $\mu (\psi (\Phi ))=0$ in ${\cal
A}_rP^{N-1}$, if $2m+1<N$. Then the closure $cl (\psi (\Phi ))$
coincides with $\psi (\Phi )\cup g(Q)$ in ${\cal A}_rP^{N-1}$. Hence
there exists $l_0\notin cl (\psi (\Phi ))$. Then consider $\pi
_{l_0}: M\to ({\cal A}_r)_{l_0}^{N-1}$. This procedure can be
prolonged, when $2m+1<N-k$, where $k$ is the number of the step of
projection. Hence $M$ can be embedded into ${\cal A}_r^{2m+1}$.
\par The approximation property follows from compactness of $M$
and the non-commutative analog of the Stone-Weierstrass theorem (see
also Theorem 2.7 in \cite{ludoyst,ludfov}).
\par Let now $M$ be a paracompact ${\cal A}_r$ holomorphic manifold
with countable atlas on $l_2(\lambda ,{\bf K})$. Spaces $l_2(\lambda
,{\cal A}_r)\oplus {\cal A}_r^m$ and $l_2(\lambda ,{\cal A}_r)\oplus
l_2(\lambda ,{\cal A}_r)$ are isomorphic as ${\cal A}_r$ Hilbert
spaces with $l_2(\lambda ,{\cal A}_r)$, since $card (\lambda )\ge
\aleph _0$. Take an additional variable $z\in {\cal A}_r$, when
$z=j\in \bf N$. Then it gives a number of a chart. Each $TU_j$ is
${\cal A}_r$ holomorphically diffeomorphic with $U_j\times
l_2(\lambda ,{\cal A}_r)$. Consider ${\cal A}_r$ holomorphic
functions $\psi $ on domains in $l_2(\lambda ,{\cal A}_r)\oplus
l_2(\lambda ,{\cal A}_r)\oplus {\cal A}_r$. Then there exists an
${\cal A}_r$ holomorphic mapping $\psi _j: M\to l_2(\lambda ,{\cal
A}_r)$ such that $\psi _j: U_j\to \psi _j(U_j)\subset l_2(\lambda
,{\cal A}_r)$ is an ${\cal A}_r$ holomorphic diffeomorphism. Then
the mapping $(\psi _1,\psi _2,...)$ provides the ${\cal A}_r$
holomorphic embedding of $M$ into $l_2(\lambda ,{\cal A}_r)$.

\end{document}